# DIFFRENCE-TYPE ESTIMATORS FOR ESTIMATION OF MEAN IN THE PRESENCE OF MEASUREMENT ERROR

[1]Viplav Kr. Singh, [†1]Rajesh Singh and [2]Florentin Smarandache

Department of Statistics, Banaras Hindu University

Varanasi-221005, India

[2]Chair of Department of Mathematics, University of New Mexico, Gallup, USA

† Corresponding author, rsinghstat@gmail.com

**Abstract**

In this paper we have suggested difference-type estimator for estimation of population mean of the study variable y in the presence of measurement error using auxiliary information. The optimum estimator in the suggested estimator has been identified along with its mean square error formula. It has been shown that the suggested estimator performs more efficient then other existing estimators. An empirical study is also carried out to illustrate the merits of proposed method over other traditional methods.

**Key Words:** Study variable, Auxiliary variable, Measurement error, Simple random Sampling, Bias, Mean Square error.

## 1. PERFORMANCE OF SUGGESTED METHOD USING SIMPLE RANDOM SAMPLING

### INTRODUCTION

The present study deals with the impact of measurement errors on estimating population mean of study variable (y) in simple random sampling using auxiliary information. In theory of survey sampling, the properties of estimators based on data are usually presupposed that the observations are the correct measurement on the characteristic being studied. When the measurement errors are negligible small, the statistical inference based on observed data continue to remain valid.

An important source of measurement error in survey data is the nature of variables (study and auxiliary). Here nature of variable signifies that the exact measurement on variables is not available. This may be due to the following three reasons:

1. The variable is clearly defined but it is hard to take correct observation at least with the currently available techniques or because of other types of practical difficulties. Eg: The level of blood sugar in a human being.

2. The variable is conceptually well defined but observation can obtain only on some closely related substitutes known as Surrogates. Eg: The measurement of economic status of a person.

3. The variable is fully comprehensible and well understood but it is not intrinsically defined. Eg: Intelligence, aggressiveness etc.

Some authors including Singh and Karpe (2008, 2009), Shalabh(1997), Allen et al. (2003), Manisha and Singh (2001, 2002), Srivastava and Shalabh (2001), Kumar et al. (2011 a,b), Malik and Singh (2013), Malik et al. (2013) have paid their attention towards the estimation of population mean $\mu_y$ of study variable using auxiliary information in the presence of measurement errors. Fuller (1995) examined the importance of measurement errors in estimating parameters in sample surveys. His major concerns are estimation of population mean or total and its standard error, quantile estimation and estimation through regression model.

**SYMBOLS AND SETUP**

Let, for a SRS scheme $(x_i, y_i)$ be the observed values instead of true values $(X_i, Y_i)$ on two characteristics (x, y), respectively for all i=(1,2,…n) and the observational or measurement errors are defined as

$$u_i = (y_i - Y_i) \tag{1}$$

$$v_i = (x_i - X_i) \tag{2}$$

where $u_i$ and $v_i$ are stochastic in nature with mean 0 and variance $\sigma_u^2$ and $\sigma_v^2$ respectively. For the sake of convenience, we assume that $u_i$'s and $v_i$'s are uncorrelated although $X_i$'s and $Y_i$'s are correlated. Such a specification can be, however, relaxed at the cost of some algebraic complexity. Also assume that finite population correction can be ignored.

Further, let the population means and variances of (x, y) be $(\mu_x, \mu_y)$ and $(\sigma_x^2, \sigma_y^2)$. $\sigma_{xy}$ and $\rho$ be the population covariance and the population correlation coefficient between x and y respectively. Also let $C_y = \frac{\sigma_y}{\mu_y}$ and $C_x = \frac{\sigma_x}{\mu_x}$ are the population coefficient of variation and $C_{yx}$ is the population coefficients of covariance in y and x.

## LARGE SAMPLE APPROXIMATION

Define:

$$e_0 = \frac{\bar{y} - \mu_y}{\mu_y} \text{ and } e_1 = \frac{\bar{x} - \mu_x}{\mu_x}$$

where, $e_0$ and $e_1$ are very small numbers and $|e_i| < 1$ $(i = 0,1)$.

Also, $E(e_i) = 0 \, (i = 0,1)$

and, $E(e_0^2) = \theta C_y^2 \left(1 + \frac{\sigma_u^2}{\sigma_y^2}\right) = \delta_0$,

$E(e_1^2) = \theta C_x^2 \left(1 + \frac{\sigma_v^2}{\sigma_x^2}\right) = \delta_1$, $E(e_0 e_1) = \theta \rho C_x C_y$, where $\theta = \frac{1}{n}$.

## 2. EXISTING ESTIMATORS AND THEIR PROPERTIES

Usual mean estimator is given by

$$\bar{y} = \sum_{i=1}^{n} \frac{y_i}{n} \tag{3}$$

Up to the first order of approximation the variance of $\bar{y}$ is given by

$$\text{Var}(\bar{y}) = \theta \mu_y^2 \left(1 + \frac{\sigma_u^2}{\sigma_y^2}\right) C_y^2 \tag{4}$$

The usual ratio estimator is given by

$$\bar{y}_R = \bar{y}\left(\frac{\mu_x}{\bar{x}}\right) \qquad (5)$$

where $\mu_x$ is known population mean of x.

The bias and MSE ($\bar{y}_R$), to the first order of approximation, are respectively, given

$$B(\bar{y}_R) = \theta \mu_y \left[\left(1 + \frac{\sigma_v^2}{\sigma_x^2}\right)C_x^2 - \rho C_y C_x\right] \qquad (6)$$

$$MSE(\bar{y}_R) = \theta \mu_y^2 \left[\left(1 + \frac{\sigma_u^2}{\sigma_y^2}\right)C_y^2 + \left(1 + \frac{\sigma_v^2}{\sigma_x^2}\right)C_x^2 - 2\rho C_y C_x\right] \qquad (7)$$

The traditional difference estimator is given by

$$\bar{y}_d = \bar{y} + k(\mu_x - \bar{x}) \qquad (8)$$

where, k is the constant whose value is to be determined.

Minimum mean square error of $\bar{y}_d$ at optimum value of

$$k = \frac{\mu_y \rho C_y}{\mu_x \left(1 + \frac{\sigma_v^2}{\sigma_x^2}\right)C_x}, \quad \text{is given by}$$

$$MSE(\bar{y}_d) = \mu_y^2 \theta \left(1 + \frac{\sigma_u^2}{\sigma_y^2}\right)C_y^2 \left[1 - \frac{\rho^2}{\left(1 + \frac{\sigma_u^2}{\sigma_y^2}\right)\left(1 + \frac{\sigma_v^2}{\sigma_x^2}\right)}\right] \qquad (9)$$

Srivastava (1967) suggested an estimator

$$\bar{y}_S = \bar{y}\left(\frac{\mu_x}{\bar{x}}\right)^{\ell_1} \qquad (10)$$

where, $\ell_1$ is an arbitrary constant.

Up to the first of approximation, the bias and minimum mean square error of $\bar{y}_S$ at optimum value of $\ell_1 = \dfrac{\rho C_y}{\left(1+\dfrac{\sigma_v^2}{\sigma_x^2}\right)C_x}$ are respectively, given by

$$B(\bar{y}_S) = \mu_y\left[\frac{\ell_1(\ell_1+1)}{2}\theta\left(1+\frac{\sigma_v^2}{\sigma_x^2}\right)C_x^2 - \ell_1\theta\rho C_y C_x\right] \tag{11}$$

$$MSE(\bar{y}_S) = \mu_y^2\,\theta\left(1+\frac{\sigma_u^2}{\sigma_y^2}\right)C_y^2\left[1 - \frac{\rho^2}{\left(1+\dfrac{\sigma_u^2}{\sigma_y^2}\right)\left(1+\dfrac{\sigma_v^2}{\sigma_x^2}\right)}\right] \tag{12}$$

Walsh (1970) suggested an estimator $\bar{y}_w$

$$\bar{y}_w = \bar{y}\left[\frac{\mu_x}{\ell_2\bar{x}+(1-\ell_2)\mu_x}\right] \tag{13}$$

where, $\ell_2$ is an arbitrary constant.

Up to the first order of approximation, the bias and minimum mean square error of $\bar{y}_w$ at optimum value of $\ell_2 = \dfrac{\rho C_y}{\left(1+\dfrac{\sigma_v^2}{\sigma_x^2}\right)C_x}$, are respectively, given by

$$B(\bar{y}_w) = \mu_y\theta\left[\ell_2^2 C_x^2\left(1+\frac{\sigma_v^2}{\sigma_x^2}\right) - \ell_2\rho C_y C_x\right] \tag{14}$$

$$MSE(\bar{y}_w) = \mu_y^2\,\theta\left(1+\frac{\sigma_u^2}{\sigma_y^2}\right)C_y^2\left[1 - \frac{\rho^2}{\left(1+\dfrac{\sigma_u^2}{\sigma_y^2}\right)\left(1+\dfrac{\sigma_v^2}{\sigma_x^2}\right)}\right] \tag{15}$$

Ray and Sahai (1979) suggested the following estimator

$$\bar{y}_{RS} = (1-\ell_3)\bar{y} + \ell_3\bar{y}\left(\frac{\bar{x}}{\mu_x}\right) \tag{16}$$

where, $\ell_3$ is an arbitrary constant.

Up to the first order of approximation, the bias and mean square of $\bar{y}_{RS}$ at optimum value of $\ell_3 = -\dfrac{\rho C_y}{\left(1+\dfrac{\sigma_v^2}{\sigma_x^2}\right)}$ are respectively, given by

$$B(\bar{y}_{RS}) = \theta\,\ell_3\,\mu_y\rho C_y C_x \tag{17}$$

$$MSE(\bar{y}_{RS}) = \mu_y^2\,\theta\left(1+\dfrac{\sigma_u^2}{\sigma_y^2}\right)C_y^2\left[1-\dfrac{\rho^2}{\left(1+\dfrac{\sigma_u^2}{\sigma_y^2}\right)\left(1+\dfrac{\sigma_v^2}{\sigma_x^2}\right)}\right] \tag{18}$$

## 3. SUGGESTED ESTIMATOR

Following Singh and Solanki (2013), we suggest the following difference-type class of estimators for estimating population mean $\overline{Y}$ of study variable y as

$$t_p = \left[\alpha_1\bar{y} + \alpha_2\bar{x}^* + (1-\alpha_1-\alpha_2)\mu_x^*\right]\left[\dfrac{\mu_x^*}{\bar{x}^*}\right]^\alpha \tag{19}$$

where $(\alpha_1,\alpha_2)$ are suitably chosen scalars such that MSE of the proposed estimator is minimum, $\bar{x}^*(=\eta\bar{x}+\lambda)$, $\mu_x^*(=\eta\mu_x+\lambda)$ with $(n,\lambda)$ are either constants or function of some known population parameters. Here it is interesting to note that some existing estimators have been shown as the members of proposed class of estimators $t_p$ for different values of $(\alpha_1,\alpha_2,\alpha,\eta,\lambda)$, which is summarized in Table 1.

**Table 1: Members of suggested class of estimators**

| Estimators | Values of Constants | | | | |
|---|---|---|---|---|---|
| | $\alpha_1$ | $\alpha_2$ | $\alpha$ | $\eta$ | $\lambda$ |
| $\bar{y}$ [Usual unbiased] | 1 | 0 | 0 | - | - |
| $\bar{y}_R$ [Usual ratio] | 1 | 0 | 1 | 1 | 0 |
| $\bar{y}_d$ [Usual difference] | 1 | $\alpha_2$ | 0 | -1 | $\mu_x$ |
| $\bar{y}_S$ [Srivastava (1967)] | 1 | 0 | $\alpha$ | 1 | 0 |
| $\bar{y}_{DS}$ [Dubey and Singh] | $\alpha_1$ | $\alpha_2$ | 0 | 1 | 0 |

The properties of suggested estimator are derived in the following theorems.

**Theorem 1.1:** Estimator $t_p$ in terms of $e_i; i = 0,1$ expressed as:

$$t_p = \left[\mu_x^* - \alpha A e_1 \mu_x^* + B\mu_x^* e_1^2 + \alpha_1 \left\{C - \alpha A C e_1 + BCe_1^2 + e_0 \mu_y - \alpha A \mu_y e_0 e_1\right\}\right.$$

$$\left. + \alpha_2 \eta \mu_x \left\{e_1 - \alpha A e_1^2\right\}\right]$$

ignoring the terms $E(e_i^r e_j^s)$ for $(r+s)>2$, where $r,s=0,1,2...$ and $i = 0,1; j = 1$ (first order of approximation).

where, $A = \dfrac{\eta \mu_x}{\eta \mu_x + \lambda}$, $B = \dfrac{\alpha(\alpha+1)}{2} A^2$ and $C = \mu_y - \mu_x^*$.

**Proof**

$$t_p = \left[\alpha_1 \bar{y} + \alpha_2 \bar{x}^* + (1 - \alpha_1 - \alpha_2)\mu_x^*\right]\left[\dfrac{\mu_x^*}{\bar{x}^*}\right]^\alpha$$

Or

$$t_p = [\alpha_1(1+e_0) + \alpha_2 \eta \mu_x e_1 + (1-\alpha_1)\mu_x^*][1 + Ae_1]^{-\alpha} \tag{20}$$

We assume $|Ae_1| < 1$, so that the term $(1+Ae_1)^{-\alpha}$ is expandable. Expanding the right hand side (20) and neglecting the terms of e's having power greater than two, we have

$$t_p = \mu_x^* - \alpha Ae_1\mu_x^* + B\mu_x^* e_1^2 + \alpha_1\{C - \alpha ACe_1 + BCe_1^2 + e_0\mu_y - \alpha A\mu_y e_0 e_1\}$$
$$+ \alpha_2\eta\mu_x\{e_1 - \alpha Ae_1^2\}$$

**Theorem: 1.2** Bias of the estimator $t_p$ is given by

$$B(t_p) = [B\mu_x^*\delta_1 + \alpha_1\{BC\delta_1 - \alpha A\mu_y\delta_{01}\} - \alpha_2\eta\mu_x A\alpha\delta_1] \tag{21}$$

**Proof:**

$$B(t_p) = E(t_p - \mu_y)$$

$$= E[\mu_x^* - \mu_y - \alpha Ae_1\mu_x^* + B\mu_x^* e_1^2 + \alpha_1\{C - \alpha ACe_1 + BCe_1^2 + e_0\mu_y - \alpha A\mu_y e_0 e_1\}$$
$$+ \alpha_2\eta\mu_x\{e_1 - \alpha Ae_1^2\}]$$

$$= [B\mu_x^*\delta_1 + \alpha_1\{BC\delta_1 - \alpha A\mu_y\delta_{01}\} - \alpha_2\eta\mu_x A\alpha\delta_1]$$

where, $\delta_0, \delta_1$ and $\delta_{01}$ are already defined in section 3.

**Theorem 1.3:** MSE of the estimator $t_p$, up to the first order of approximation is

$$MSE(t_p) = \alpha_1^2\{C^2 + \mu_y^2\delta_0 + \delta_1(\alpha^2 A^2 C^2 + 2BC^2) - 4\alpha AC\mu_y\delta_{01}\} + \alpha_2^2\eta^2\mu_x^2\delta_1$$
$$+ \{C^2 + \delta_1(\alpha^2 A^2\mu_x^2 - 2BC\mu_x^*)\} - 2\alpha_1\{C^2 + \delta_1(BC^2 - BC\mu_x^* - \alpha^2 A^2 C\mu_x^*) + \delta_{01}\alpha A\mu_y(\mu_x^* - C)\}$$
$$- 2\alpha_2\eta\mu_x\alpha A\delta_1(\mu_x^* - C) + 2\alpha_1\alpha_2\eta\mu_x(\mu_y\delta_{01} - 2A\alpha C\delta_1) \tag{22}$$

**Proof:**

$$MSE(t_p) = E(t_p - \mu_y)^2$$

$$= E[\alpha_1\{C - A\alpha Ce_1 + e_0\mu_y + BCe_1^2 - \alpha A\mu_y e_0 e_1\} + \alpha_2\eta\mu_x\{e_1 - A\alpha e_1^2\}$$

$$- C + \alpha Ae_1\mu_x^* - B\mu_x^* e_1^2]^2$$

Squaring and then taking expectations of both sides, we get the MSE of the suggested estimator up to the first order of approximation as

$$MSE(t_p) = \alpha_1^2\{C^2 + \mu_y^2\delta_0 + \delta_1(\alpha^2A^2C^2 + 2BC^2) - 4\alpha AC\mu_y\delta_{01}\} + \alpha_2^2\eta^2\mu_x^2\delta_1$$

$$+ \{C^2 + \delta_1(\alpha^2A^2\mu_x^2 - 2BC\mu_x^*)\} - 2\alpha_1\{C^2 + \delta_1(BC^2 - BC\mu_x^* - \alpha^2A^2C\mu_x^*) + \delta_{01}\alpha A\mu_y(\mu_x^* - C)\}$$

$$- 2\alpha_2\eta\mu_x\alpha A\delta_1(\mu_x^* - C) + 2\alpha_1\alpha_2\eta\mu_x(\mu_y\delta_{01} - 2A\alpha C\delta_1)$$

Equation (22) can be written as:

$$MSE(t_p) = \alpha_1^2\varphi_1 + \alpha_2^2\varphi_2 - 2\alpha_1\varphi_3 - 2\alpha_2\varphi_4 + 2\alpha_1\alpha_2\varphi_5 + \varphi \tag{23}$$

Differentiating (23) with respect to $(\alpha_1, \alpha_2)$ and equating them to zero, we get the optimum values of $(\alpha_1, \alpha_2)$ as

$$\alpha_{1(opt)} = \frac{\varphi_2\varphi_3 - \varphi_4\varphi_5}{\varphi_1\varphi_2 - \varphi_5^2} \text{ and } \alpha_{2(opt)} = \frac{\varphi_1\varphi_4 - \varphi_3\varphi_5}{\varphi_1\varphi_2 - \varphi_5^2}$$

where, $\varphi_1 = C^2 + \mu_y^2\delta_0 + \delta_1(\alpha^2A^2C^2 + 2BC^2) - 4\alpha AC\mu_y\delta_{01}$

$$\varphi_2 = \eta^2\mu_x^2\delta_1$$

$$\varphi_3 = C^2 + \delta_1(BC^2 - BC\mu_x^* - \alpha^2A^2C\mu_x^*) + \delta_{01}\alpha A\mu_y(\mu_x^* - C)$$

$$\varphi_4 = \eta\mu_x\alpha A\delta_1(\mu_x^* - C)$$

$$\varphi_5 = \eta\mu_x(\mu_y\delta_{01} - 2A\alpha C\delta_1)$$

$$\varphi = C^2 + \delta_1(\alpha^2A^2\mu_x^2 - 2BC\mu_x^*)$$

In the Table 2 some estimators are listed which are particular members of the suggested class of estimators $t_p$ for different values of $(\alpha, \eta, \lambda)$.

**Table 2:** Particular members of the suggested class of estimators $t_p$

| Estimators | Values of constants | | |
|---|---|---|---|
| | $\alpha$ | $\eta$ | $\lambda$ |
| $t_1 = [\alpha_1 \bar{y} + \alpha_2 \bar{x} + (1-\alpha_1 - \alpha_2)\mu_x] \left[ \dfrac{\mu_x}{\bar{x}} \right]$ | -1 | 1 | 0 |
| $t_2 = [\alpha_1 \bar{y} + \alpha_2 (\bar{x}+1) + (1-\alpha_1 - \alpha_2)(\mu_x+1)] \left[ \dfrac{\mu_x+1}{\bar{x}+1} \right]$ | 1 | 1 | 1 |
| $t_3 = [\alpha_1 \bar{y} + \alpha_2 (\bar{x}+1) + (1-\alpha_1 - \alpha_2)(\mu_x+1)] \left[ \dfrac{\mu_x}{\bar{x}} \right]^{-1}$ | -1 | 1 | 1 |
| $t_4 = [\alpha_1 \bar{y} + \alpha_2 (\bar{x}+\rho) + (1-\alpha_1 - \alpha_2)(\mu_x+\rho)] \left[ \dfrac{\mu_x+\rho}{\bar{x}+\rho} \right]^{-1}$ | -1 | 1 | $\rho$ |
| $t_5 = [\alpha_1 \bar{y} + \alpha_2 (\bar{x}+C_x) + (1-\alpha_1 - \alpha_2)(\mu_x+C_x)] \left[ \dfrac{\mu_x+C_x}{\bar{x}+C_x} \right]^{-1}$ | -1 | 1 | $C_x$ |
| $t_6 = [\alpha_1 \bar{y} + \alpha_2 (\bar{x}-C_x) + (1-\alpha_1 - \alpha_2)(\mu_x-C_x)] \left[ \dfrac{\mu_x-C_x}{\bar{x}-C_x} \right]$ | -1 | 1 | $-C_x$ |
| $t_7 = [\alpha_1 \bar{y} - \alpha_2 (\bar{x}+1) - (1-\alpha_1 - \alpha_2)(\mu_x+C_x)] \left[ \dfrac{\mu_x+C_x}{\bar{x}+C_x} \right]^{-1}$ | -1 | -1 | -1 |

## 4. EMPIRICAL STUDY

**Data statistics:** The data used for empirical study has been taken from Gujarati (2007)

Where, $Y_i$ =True consumption expenditure,

$X_i$ =True income,

$y_i$ =Measured consumption expenditure,

$x_i$ = Measured income.

| n | $\mu_y$ | $\mu_x$ | $\sigma_y^2$ | $\sigma_x^2$ | $\rho$ | $\sigma_u^2$ | $\sigma_v^2$ |
|---|---|---|---|---|---|---|---|
| 10 | 127 | 170 | 1278 | 3300 | 0.964 | 36 | 36 |

The percentage relative efficiencies (PRE) of various estimators with respect to the mean per unit estimator of $\bar{Y}$, that is $\bar{y}$, can be obtained as

$$PRE(.) = \frac{Var(\bar{y})}{MSE(.)} * 100$$

**Table 3: MSE and PRE of estimators with respect to $\bar{y}$**

| Estimators | Mean Square Error | Percent Relative Efficiency |
|---|---|---|
| $\bar{y}$ | 131.4 | 100 |
| $\bar{y}_R$ | 21.7906 | 603.0118 |
| $\bar{y}_d$ | 13.916 | 944.1285 |
| $\bar{y}_S$ | 13.916 | 944.1285 |
| $\bar{y}_{DS}$ | 13.916 | 944.1285 |
| $t_1$ | 10.0625 | 1236.648 |
| $t_2$ | 9.92677 | 1323.693 |
| $t_3$ | **6.82471** | 1925.356 |
| $t_4$ | 6.9604 | 1887.818 |
| $t_5$ | 9.3338 | 1407.774 |
| $t_6$ | 11.9246 | 1101.923 |
| $t_7$ | 7.9917 | 1644.194 |

## 5. PERFORMANCE OF SUGGESTED ESTIMATOR IN STRATIFIED RANDOM SAMPLING

### SYMBOLS AND SETUP

Consider a finite population $U = (u_1, u_2, ..., u_N)$ of size N and let X and Y respectively be the auxiliary and study variables associated with each unit $u_j = (j = 1,2,......, N)$ of population. Let the population of N be stratified in to L strata with the $h^{th}$ stratum containing $N_h$ units, where h=,1,2,3,....,L such that $\sum_{i=1}^{L} N_h = N$. A simple random size $n^h$ is drown without replacement from the $h^{th}$ stratum such that $\sum_{i=1}^{L} n_h = n$. Let $(y_{hi}, X_{hi})$ of two characteristics (Y,X) on $i^{th}$ unit of the $h^{th}$ stratum, where i=,1,2,...,$N_h$. In addition let

$$(\bar{y}_h = \frac{1}{n_h} \sum_{i=1}^{n_h} y_{hi}, \bar{x}_h = \frac{1}{n_h} \sum_{i=1}^{n_h} x_{hi}),$$

$$(\bar{y}_{st} = \sum_{i=1}^{n_h} W_h \bar{y}_h, \bar{x}_{st} = \sum_{i=1}^{n_h} W_h \bar{x}_h),$$

$$(\mu_{Yh} = \frac{1}{N_h} \sum_{i=1}^{N_h} y_{hi}, \mu_{Xh} = \frac{1}{N_h} \sum_{i=1}^{N_h} x_{hi}),$$

And $(\mu_Y = \sum_{h=1}^{L} W_h \mu_{Yh}, \mu_X = \sum_{i=1}^{N_h} W_h \mu_{Xh})$ be the samples means and population means of (Y,X) respectively, where $W_h = \frac{N_h}{N}$ is the stratum weight. Let the observational or measurement errors be

$$u_{hi} = y_{hi} - Y_{hi} \quad (24)$$

$$v_{hi} = x_{hi} - X_{hi} \quad (25)$$

Where $u_{hi}$ and $v_{hi}$ are stochastic in nature and are uncorrelated with mean zero and variances $\sigma_{Vh}^2$ and $\sigma_{Uh}^2$ respectively. Further let $\rho_h$ be the population correlation coefficient between Y and X in the $h^{th}$ stratum. It is also assumed that the finite population correction terms $(1-f_h)$ and $(1-f)$ can be ignored where $f_h = \frac{n_h}{N_h}$ and $f = \frac{n}{N}$.

## LARGE SAMPLE APPROXIMATION

Let,

$$\bar{y}_{st} = \mu_Y(1+e_{0h}), \text{ and } \bar{x}_{st} = \mu_X(1+e_{1h})$$

such that, $E(e_{0h}) = E(e_{1h}) = 0$,

$$E(e_{0h}^2) = \frac{C_{Yh}^2}{n_h}\left(1+\frac{\sigma_{Uh}^2}{\sigma_{Yh}^2}\right) = \frac{C_{Yh}^2}{n_h \theta_{Yh}} = \nabla_0,$$

$$E(e_{1h}^2) = \frac{C_{Xh}^2}{n_h}\left(1+\frac{\sigma_{Vh}^2}{\sigma_{Xh}^2}\right) = \frac{C_{Xh}^2}{n_h \theta_{Xh}} = \nabla_1,$$

$$E(e_{0h} e_{1h}) = \frac{1}{n_h}\rho_h C_{Yh} C_{Xh} = \nabla_{01}.$$

where, $C_{Yh} = \frac{\sigma_{Yh}}{\mu_{Yh}}$, $C_{Xh} = \frac{\sigma_{Xh}}{\mu_{Xh}}$, $\theta_{Yh} = \frac{\sigma_{Uh}^2}{\sigma_{Uh}^2 + \sigma_{Yh}^2}$ and $\theta_{Xh} = \frac{\sigma_{Vh}^2}{\sigma_{Vh}^2 + \sigma_{Xh}^2}$.

## EXISTING ESTIMATORS AND THEIR PROPERTIES

$\bar{y}_{st}$ is usual unbiased estimator in stratified random sampling scheme.

The usual combined ratio estimator in stratified random sampling in the presence of measurement error is defined as-

$$T_R = \bar{y}_{st} \frac{\mu_x}{\bar{x}_{st}} \tag{26}$$

The usual combined product estimator in the presence of measurement error is defined as-

$$T_{PR} = \bar{y}_{st} \frac{\bar{x}_{st}}{\mu_x} \tag{27}$$

Combined difference estimator in stratified random sampling is defined in the presence of measurement errors for a population mean, as

$$T_D = \bar{y}_{st} + d(\mu_x - \bar{x}_{st}) \tag{28}$$

The variance and mean square term of above estimators, up to the first order of approximation, are respectively given by

$$\text{Var}(\bar{y}_{st}) = \frac{C_{Xh}^2}{n_h}\left(1 + \frac{\sigma_{Uh}^2}{\sigma_{Yh}^2}\right) \tag{29}$$

$$\text{MSE}(T_R) = \sum_{h=1}^{L} \frac{W_h^2}{n_h}\left[\frac{\sigma_{Yh}^2}{\theta_{Yh}} + R\left(\frac{\sigma_{Xh}^2}{\theta_{Xh}}\right)\right](R - 2\beta_{YXh}\theta_{Xh}) \tag{30}$$

$$\text{MSE}(T_P) = \sum_{h=1}^{L} \frac{W_h^2}{n_h}\left[\frac{\sigma_{Yh}^2}{\theta_{Yh}} + R\left(\frac{\sigma_{Xh}^2}{\theta_{Xh}}\right)\right](R + 2\beta_{YXh}\theta_{Xh}) \tag{31}$$

$$\text{MSE}(T_D) = \sum_{h=1}^{L} \frac{W_h^2}{n_h}\left(\frac{\sigma_{Yh}^2}{\theta_{Yh}}\right) + d^2 \sum_{h=1}^{L} \frac{W_h^2}{n_h}\left(\frac{\sigma_{Xh}^2}{\theta_{Xh}}\right) - 2d\sum_{h=1}^{L} \frac{W_h^2}{n_h}\beta_{XYh}\sigma_{Xh}^2 \tag{32}$$

$$\text{where, } d_{opt} = \frac{\sum_{h=1}^{L} \frac{W_h^2}{n_h}\beta_{XYh}\sigma_{Xh}^2}{\sum_{h=1}^{L} \frac{W_h^2}{n_h}\left(\frac{\sigma_{Xh}^2}{\theta_{Xh}}\right)}$$

## 6. SUUGESTED ESTIMATOR AND ITS PROPERTIES

Let $B(.)$ and $M(.)$ denote the bias and mean square error (M.S.E) of an estimator under given sampling design. Estimator $t_p$ defined in equation (19) can be written in stratified random sampling as

$$T_P = \left[\beta_1 \bar{y}_{st} + \beta_2 \bar{x}_{st}^* + (1 - \beta_1 - \beta_2)\mu_x^*\right]\left[\frac{\mu_x^*}{\bar{x}_{st}^*}\right]^\beta \tag{33}$$

where $(\alpha_1, \alpha_2)$ are suitably chosen scalars such that MSE of proposed estimator is minimum, $\bar{x}_{st}^*(= \eta\bar{x}_{st} + \lambda), \mu_x^*(= \eta\mu_x + \lambda)$ with $(n, \lambda)$ are either constants or functions of some known population parameters. Here it is interesting to note that some existing estimators have been found particular members of proposed class of estimators $T_P$ for different values of $(\alpha_1, \alpha_2, \alpha, \eta, \lambda)$, which are summarized in Table 4.

**Table 4: Members of proposed class of estimators $T_P$**

| Estimators | Values of Constants | | | | |
|---|---|---|---|---|---|
| | $\alpha_1$ | $\alpha_2$ | $\alpha$ | $\eta$ | $\lambda$ |
| $\bar{y}_{st}$ [Usual unbiased] | 1 | 0 | 0 | - | - |
| $T_R$ [Usual ratio] | 1 | 0 | 1 | 1 | 0 |
| $T_{PR}$ [Usual product] | 1 | 0 | -1 | 1 | 0 |
| $T_D$ [Usual difference] | 1 | $\alpha_2$ | 0 | -1 | $\mu_x$ |

**Theorem 2.1:** Estimator $T_P$ in terms of $e_i; i = 0,1$ by ignoring the terms $E(e_{ih}^r e_{jh}^s)$ for $(r+s)>2$, where $r,s=0,1,2...$ and $i = 0,1; j = 1$, can be written as

$$T_P = [\mu_x^* - \beta A' e_{1h} \mu_x^* + B' \mu_x^* e_{1h}^2 + \beta_1 \{C'' - \beta A' C e_{1h} + B' C' e_{1h}^2 + e_{0h} \mu_y - \beta A' \mu_y e_{0h} e_{1h}\}$$

$$+ \beta_2 \eta \mu_x \{e_{1h} - \beta A' e_{1h}^2\}]$$

where, $A' = \dfrac{\eta \mu_x}{\eta \mu_x + \lambda}$, $B' = \dfrac{\beta(\beta+1)}{2} A'^2$ and $C' = \mu_y - \mu_x^*$.

**Proof**

$$T_P = \left[\beta_1 \bar{y}_{st} + \beta_2 \bar{x}_{st}^* + (1-\beta_1-\beta_2)\mu_x^*\right]\left[\dfrac{\mu_x^*}{\bar{x}_{st}^*}\right]^\beta$$

$$= [\beta_1(1+e_{0h}) + \beta_2 \eta \mu_x e_{1h} + (1-\beta_1)\mu_x^*][1+A'e_{1h}]^{-\beta} \qquad (34)$$

We assume $|A'e_{1h}| < 1$, so that the term $(1+A'e_{1h})^{-\beta}$ is expandable. Thus by expanding the right hand side (20) and neglecting the terms of e's having power greater than two, we have

$$T_P = [\mu_x^* - \beta A'e_{1h}\mu_x^* + B'\mu_x^* e_{1h}^2 + \beta_1\{C' - \beta A'C'e_{1h} + B'C'e_{1h}^2 + e_{0h}\mu_y - \beta A'\mu_y e_{0h} e_{1h}\}$$
$$+ \beta_2\eta\mu_x\{e_{1h} - \beta A'e_{1h}^2\}]$$

**Theorem: 2.2** Bias of $T_p$ is given by

$$B(T_P) = [B'\mu_x^*\nabla_1 + \beta_1\{B'C'\nabla_1 - \beta A'\mu_y\nabla_{01}\} - \beta_2\eta\mu_x A'\beta\nabla_1] \tag{35}$$

**Proof:**

$$B(T_P) = E(T_P - \mu_y)$$

$$= E[\mu_x^* - \mu_y - \beta A'e_{1h}\mu_x^* + B'\mu_x^* e_{1h}^2 + \beta_1\{C' - \beta A'C'e_{1h} + B'C'e_{1h}^2 + e_{0h}\mu_y - \beta A'\mu_y e_{0h} e_{1h}\}$$
$$+ \beta_2\eta\mu_x\{e_{1h} - \beta A'e_{1h}^2\}]$$

$$= [B'\mu_x^*\nabla_1 + \beta_1\{B'C'\nabla_1 - \beta A'\mu_y\nabla_{01}\} - \beta_2\beta\eta\mu_x A'\nabla_1]$$

where, $\nabla_0, \nabla_1$ and $\nabla_{01}$ are already defined in section 3.

**Theorem: 2.3** Mean square error of $T_P$, up to the first order of approximation is given by

$$MSE(T_P) = \beta_1^2\{C'^2 + \mu_y^2\nabla_0 + \nabla_1(\beta^2 A'^2 C'^2 + 2B'C'^2) - 4\beta A'C'\mu_y\nabla_{01}\} + \beta_2^2\eta^2\mu_x^2\nabla_1$$

$$+ \{C'^2 + \nabla_1(\beta^2 A'^2 \mu_x^2 - 2B'C'\mu_x^*)\} - 2\beta_1\{C'^2 + \nabla_1(B'C'^2 - B'C'\mu_x^* - \beta^2 A'^2 C'\mu_x^*) + \nabla_{01}\beta A\mu_y(\mu_x^* - C')\}$$
$$- 2\beta_2\eta\mu_x\beta A\nabla_1(\mu_x^* - C') + 2\beta_1\beta_2\eta\mu_x(\mu_y\nabla_{01} - 2A'\beta C'\nabla_1) \tag{36}$$

**Proof:**

$$MSE(T_P) = E(T_P - \mu_y)^2$$

$$MSE(T_P) = \beta_1^2\{C'^2 + \mu_y^2\nabla_0 + \nabla_1(\beta^2 A'^2 C'^2 + 2B'C'^2) - 4\beta A'C'\mu_y\nabla_{01}\} + \beta_2^2\eta^2\mu_x^2\nabla_1$$

$$+ \{C'^2 + \nabla_1(\beta^2 A'^2 \mu_x^2 - 2B'C'\mu_x^*)\} - 2\beta_1\{C'^2 + \nabla_1(B'C'^2 - B'C'\mu_x^* - \beta^2 A'^2 C'\mu_x^*) + \nabla_{01}\beta A\mu_y(\mu_x^* - C')\}$$
$$- 2\beta_2\eta\mu_x\beta A\nabla_1(\mu_x^* - C') + 2\beta_1\beta_2\eta\mu_x(\mu_y\nabla_{01} - 2A'\beta C'\nabla_1)$$

MSE($T_p$) can also be written as

$$\text{MSE}(T_P) = \beta_1^2 \chi_1 + \beta_2^2 \chi_2 - 2\beta_1\chi_3 - 2\beta_2\chi_4 + 2\beta_1\beta_2\chi_5 + \chi \tag{37}$$

Differentiating equation (37) with respect to $(\beta_1, \beta_2)$ and equating it to zero, we get the optimum values of $(\beta_1, \beta_2)$ respectively, as

$$\beta_{1(opt)} = \frac{\chi_2\chi_3 - \chi_4\chi_5}{\chi_1\chi_2 - \chi_5^2} \text{ and } \beta_{2(opt)} = \frac{\chi_1\chi_4 - \chi_3\chi_5}{\chi_1\chi_2 - \chi_5^2}$$

where,
$$\chi_1 = C'^2 + \mu_y^2 \nabla_0 + \nabla_1(\beta^2 A'^2 C'^2 + 2B'C'^2) - 4\beta A'C'\mu_y\nabla_{01}$$

$$\chi_2 = \eta^2 \mu_x^2 \nabla_1$$

$$\chi_3 = C'^2 + \nabla_1(BC^2 - B'C'\mu_x^* - \beta^2 A'^2 C'\mu_x^*) + \nabla_{01}\beta A'\mu_y(\mu_x^* - C')$$

$$\chi_4 = \eta\mu_x\beta A'\nabla_1(\mu_x^* - C')$$

$$\chi_5 = \eta\mu_x(\mu_y\nabla_{01} - 2A'\beta C'\nabla_1)$$

$$\chi = C'^2 + \nabla_1(\beta^2 A'^2 \mu_x^2 - 2B'C'\mu_x^*)$$

With the help of these values, we get the minimum MSE of the suggested estimator $T_p$.

### 7. DISCUSSION AND CONCLUSION

In the present study, we have proposed difference-type class of estimators of the population mean of a study variable when information on an auxiliary variable is known in advance. The asymptotic bias and mean square error formulae of suggested class of estimators have been obtained. The asymptotic optimum estimator in the suggested class has been identified with its properties. We have also studied some traditional methods of estimation of population mean in the presence of measurement error such as usual unbiased, ratio, usual difference estimators suggested by Srivastava(1967), dubey and singh( 2001), which are found to be particular members of suggested class of estimators. In addition, some new members of suggested class of estimators have also been generated in simple random sampling case. An empirical study is carried to throw light on the performance of suggested estimators over other existing estimators using simple random sampling scheme. From the Table 3, we observe that suggested estimator $t_3$ performs better than the other estimators considered in the present study and which reflects the usefulness of suggested method in practice.

**REFERENCES**

Allen, J., Singh, H. P. and Smarandache, F. (2003): A family of estimators of population mean using mutliauxiliary information in presence of measurement errors. International Journal of Social Economics 30(7), 837–849.

A.K. Srivastava and Shalabh (2001). Effect of measurement errors on the regression method of estimation in survey sampling. *Journal of Statistical Research*, Vol. **35**, No. 2, pp. 35-44.

Bahl, S. and Tuteja, R. K. (1991): Ratio and product type exponential estimator. Information and optimization sciences12 (1), 159-163.

Chandhok, P.K., & Han, C.P.(1990):On the efficiency of the ratio estimator under Midzuno scheme with measurement errors. Journal of the Indian statistical Association,28,31-39.

Dubey, V. and Singh, S.K. (2001). An improved regression estimator for estimating population mean, J. Ind. Soc. Agri. Statist., 54, p. 179-183.

Gujarati, D. N. and Sangeetha (2007): Basic econometrics. Tata McGraw – Hill.

Koyuncu, N. and Kadilar, C. (2010): On the family of estimators of population mean in stratified sampling. Pakistan Journal of Statistics. Pak. J. Stat. 2010 vol 26 (2), 427-443.

Kumar, M. , Singh, R., Singh, A.K. and Smarandache, F. (2011 a): Some ratio type estimators under measurement errors. WASJ 14(2) :272-276.

Kumar, M., Singh, R., Sawan, N. and Chauhan, P. (2011b): Exponential ratio method of estimators in the presence of measurement errors. Int. J. Agricult. Stat. Sci. 7(2): 457-461.

Malik, S. and Singh, R. (2013) : An improved class of exponential ratio- type estimator in the presence of measurement errors. OCTOGON Mathematical Magazine, 21,1, 50-58.

Malik, S., Singh, J. and Singh, R. (2013) : A family of estimators for estimating the population mean in simple random sampling under measurement errors. JRSA, 2(1), 94-101.

Manisha and Singh, R. K. (2001): An estimation of population mean in the presence of measurement errors. Journal of Indian Society of Agricultural Statistics 54(1), 13–18.

Manisha and Singh, R. K. (2002): Role of regression estimator involving measurement errors. Brazilian journal of probability Statistics 16, 39- 46.

Shalabh (1997): Ratio method of estimation in the presence of measurement errors. Journal of Indian Society of Agricultural Statistics 50(2):150– 155.

Singh, H. P. and Karpe, N. (2008): Ratio-product estimator for population mean in presence of measurement errors. Journal of Applied Statistical Sciences 16, 49–64.

Singh, H. P. and Karpe, N. (2009): On the estimation of ratio and product of two populations means using supplementary information in presence of measurement errors. Department of Statistics, University of Bologna, 69(1), 27-47.

Singh, H. P. and Vishwakarma, G. K. (2005): Combined Ratio-Product Estimator of Finite Population Mean in Stratified Sampling. Metodologia de Encuestas 8: 35- 44.

Singh, H. P., Rathour, A., Solanki, R. S.(2013): An improvement over difference method of estimation of population mean. JRSS, 6(1):35-46.

Singh H. P., Rathour A., Solanki R.S. (2013): an improvement over difference method of estimation of population mean. JRSS, 6(1): 35-46.

Srivastava, S. K. (1967): An estimator using auxiliary information in the sample surveys. Calcutta statistical Association Bulletin 16,121-132.

Walsh, J.E.(1970). Generalisation of ratio estimate for population total. Sankhya A.32, 99-106.